\documentclass{amsart}
\usepackage[latin1]{inputenc}
\usepackage{amsfonts,amssymb,amsmath,formel}
\newcommand{\ETA}{{\boldmath\mbox{$\eta$}}}
\newcommand{\SIGMA}{{\boldmath\mbox{$\sigma$}}}

\begin{document}
\begin{center}
{\bf Random diophantine equations of additive type}
\\[3pt]
J. Br\"udern and R. Dietmann
\end{center}

\begin{center}
{\bf I. Introduction.}
\end{center}

\setcounter{section}{1}\setcounter{equation}{0}
 In this memoir, we investigate the solubility
of diagonal diophantine equations
\be{1.1}
a_1 x_1^k + a_2 x_2^k + \ldots + a_s x_s^k =0,
\ee
and the distribution of their solutions. This is a theme that has
received much interest in the past (see Vaughan \cite{hlm}, Vaughan
and Wooley \cite{survey}, Heath-Brown 
\cite{HB},
Swinnerton-Dyer \cite{SD} and the extensive bibliographies in
\cite{hlm,survey}). Our main concern is with the validity of
the Hasse principle, and with a bound for the smallest non-zero
solution in integers whenever such a solution exists. The approach is
of a statistical nature. Very roughly speaking, we shall show that
whenever $s>4k$ and the vector ${\bf a}= (a_1, \ldots, a_s)\in\ZZ^s$
is chosen at random, then almost surely the Hasse principle holds for
\rf{1.1}, and if there are solutions in integers, not all zero, then
there is one with $|{\bf x}|\ll |{\bf a}|^{2/(s-2k-2)}$. Here and
later, we write  $|{\bf
  x}|=\max |x_j|$.
\medskip

We now set the scene to describe our results in precise form. To avoid
trivialities, suppose throughout that $k\in\NN$, $k\ge 2$, and that
$a_j\in\ZZ\bs\{0\}$. Since \rf{1.1} has the trivial solution ${\bf
  x}={\bf 0}$, it will be convenient to describe the equation \rf{1.1}
as {\it soluble} over a given field if there exists a solution in that
field other than the trivial one. If \rf{1.1} is soluble over $\RR$
and over $\QQ_p$
for all primes $p$,  then \rf{1.1} is called {\it
  locally soluble}. We denote by ${\cal C}={\cal C}(k,s)$ the set
of all ${\bf a}$ with $a_j\in\ZZ\bs \{0\}$ for which \rf{1.1} is
locally soluble. Note that whenever \rf{1.1} is soluble over $\QQ$,
then ${\bf a}\in{\cal C}$. The inverse conclusion is known as the {\it
  Hasse principle} for the equation \rf{1.1}. Recall that when $k=2$,
then the Hasse principle holds for any $s$, as a special case of the
Hasse-Minkowski theorem.

Whenever $s>2k$, a formal use of the Hardy-Littlewood method leads one
to expect an asymptotic formula for the number $\rho_{\bf a}(B)$ of
solutions of \rf{1.1} in integers $x_j$ within the box $|{\bf x}|\le
B$. This takes the shape
\be{1.2}
\rho_{\bf a} (B) = B^{s-k} J_{\bf a} \prod_{p} \chi_p({\bf
  a}) + o(B^{s-k}) \quad (B\to\inf)
\ee
where
\be{1.3}
\chi_p ({\bf a}) = \lim_{h\to\inf} p^{h(1-s)} \#\{1\le x_j \le p^h:
a_1 x_1^k + \ldots + a_s x_s^k \equiv 0 \bmod p^h\}
\ee
is a measure for the density of the solutions of \rf{1.1} in $\QQ_p$,
and similarly, $ J_{\bf a} $ is related to the surface area of the
real solutions of \rf{1.1} within the box $[-1,1]^s$. A precise
definition of $J_{\bf a} $ is given in \rf{3.7} below.  

As we shall see later, a condition milder than the 
current hypothesis $s>2k$ suffices to
confirm that the limits \rf{1.3} exist for all primes $p$, and that
the Euler product
\be{1.4}
{\fk S}_{\bf a}= \prod_p \chi_p({\bf a})
\ee
is absolutely convergent. Moreover, an
application of Hensel's lemma shows that $\chi_p({\bf a})$ is positive
if and only if \rf{1.1} is soluble in $\QQ_p$. Likewise, one finds
that $J_{\bf a}$ is positive if and only if \rf{1.1} is
soluble over $\RR$. It follows that \rf{1.1} is locally soluble if and
only if
\be{1.5}
J_{\bf a} {\fk S}_{\bf a}>0.
\ee
Consequently, if \rf{1.2} holds, then the equation \rf{1.1} obeys the
Hasse principle.
\medskip

The validity of \rf{1.2}, and hence of  the Hasse principle for the
underlying diophantine equations, is regarded to be a save conjecture
in the range $s>2k$, and in the special case $k=2$, $s>4$ rigorous
proofs of \rf{1.2} are 
available by various methods (see chapter 2 of \cite{hlm}
for one approach). When $k=3$, the formula \rf{1.2} is known to hold
whenever $s\ge 8$ (implicit in Vaughan \cite{V86}), and the Hasse
principle holds for $s\ge 7$ (Baker \cite{B89}). For larger $k$, much
less is known. The asymptotic formula \rf{1.2} has been established
when $s\ge \f{1}{2} k^2\log k (1+o(1))$, and the Hasse principle may
be verified when $s\ge k \log k (1+o(1))$, see Ford \cite{F} and
Wooley \cite{W92}. Although these results fall short of the expected
one by a factor of $\log k$ at least, with respect to the number of
variables, it seems difficult to establish \rf{1.2} on average over
${\bf a}$ when $s$ is significantly smaller than in the aforementioned
work of Ford. However, one may choose $B$ as a suitable
function of $|{\bf a}|$, say $B=|{\bf a}|^\theta$, and then investigate
whether \rf{1.5} holds for almost all ${\bf a}$. This approach is
successful whenever $s>4k$ and $\theta$ is approximately as large as
$2/(s-2k)$, and suffices to confirm the conclusions alluded to in
the introductory paragraph. The principal step is contained in the
following mean value theorem. Before this is formulated, recall that 
${\bf a}$ is reserved for integral vectors with non-zero
entries; this convention applies within the summation below, and
elsewhere in this paper. Also, when $s$ is a natural number, let $\hat
s$ denote the largest even
  integer strictly smaller than $s$.
\\[2ex]
{\sc Theorem 1.1.} {\it Let $k\ge 3$ and $s>4k$. Then there is a 
positive number
  $\del$ such that whenever $A,B$ are real numbers satisfying
\be{1.6}
1\le 2B^k \le A \le B^{({\hat s}-2k)/2},
\ee
one has}
$$
\sum_{|{\bf a}|\le A} |\rho_{\bf a} (B) - J_{\bf a}
{\fk S}_{\bf a} B^{s-k}| \ll A^{s-1-\del} B^{s-k}.
$$

This theorem actually remains valid when $k=2$, but the proof we give below
needs some adjustments. We have excluded $k=2$ from the discussion mainly 
because in that particular case one can say more, by different 
methods. Hence,
from now on, we assume throughout that $k\ge 3$.

As a simple corollary, we note that subject to the conditions in
Theorem 1.1, the number of ${\bf a}$ with $|{\bf a}|\le A$ for which
the inequality
\be{1.7}
|\rho_{\bf a}(B) -J_{\bf a}
{\fk S}_{\bf a} B^{s-k}  |>
|{\bf a}|^{-1} B^{s-k-\del}
\ee
holds, does not exceed $O(A^{s-\f{1}{2}\del})$. To deduce the Hasse
principle for those ${\bf a}$ where \rf{1.7} fails, one needs a lower
bound for $J_{\bf a} {\fk S}_{\bf a}$ whenever this number
is non-zero. When $k$ is odd, \rf{1.1} is soluble over $\RR$, and one
may show that
\be{1.8}
J_{\bf a} \gg |{\bf a}|^{-1}
\ee
holds for all ${\bf a}$. When $k$ is even, \rf{1.1} is soluble over
$\RR$ if and only if the $a_j$ are not all of the same sign, and if
this is the case, then again \rf{1.8} holds. These facts will be
demonstrated in \S3. For the ``singular product'' we have the following
result.
\\[2ex]
{\sc Theorem 1.2.} {\it Let $s\ge k+3$, and let $\eta$ be a positive
  number. Then there exists a positive number $\gm$ such that}
$$
\# \{|{\bf a}| \le A: \, 0 < {\fk S}_{\bf a} < A^{-\eta}\} \ll
A^{s-\gm}.
$$

\medskip
We are ready to derive the main result. Let $s> 4k$, and let $\del$ be
the positive number supplied by Theorem 1.1. Suppose that ${\bf a}
\in{\cal C}(k,s)$ satisfies $\f{1}{2} A < |{\bf a}| \le A$, and choose
$B = A^{2/({\hat s}-2k)}$ in accordance with \rf{1.6}. In Theorem 1.2, we
take $\eta = \del / ({\hat s}-2k)$ so that $A^\eta = B^{\del/2}$. If ${\bf
  a}$ is not counted in Theorem 1.2, then ${\fk S}_{\bf a} \ge
A^{-\eta}$, and if ${\bf a}$ also violates \rf{1.7}, then by \rf{1.8}
one has 
$$
\rho_{\bf a} (B) \ge J_{\bf a}
{\fk S}_{\bf a} B^{s-k} -
|{\bf a}|^{-1} B^{s-k-\del} \gg B^{s-k} A^{-1-\eta} \gg A^{1-\eta}.
$$
It follows that \rf{1.1} has an integral solution with $0<|{\bf x}|\le
B \ll |{\bf a}|^{2/({\hat s}-2k)}$, for these choices of ${\bf a}$. The remaining
${\bf a}\in{\cal C}(k,s)$ with $\f{1}{2}A < |{\bf a}|\le A$ are counted
in \rf{1.7} or in Theorem 1.2. Therefore, there are at most
$O(A^{s-\min (\del, \gm)})$ such ${\bf a}$. We now sum for $A$ over
powers of $2$ to conclude as follows.
\\[2ex]
{\sc Theorem 1.3.} {\it Let $s > 4k$. Then, there is a positive number
  $\theta$ such that the number of ${\bf a}\in
  {\cal C}(k,s)$ for which the equation \rf{1.1} has no integral
  solution in the range $0 < |{\bf x}| \le |{\bf a}|^{2/({\hat s}-2k)}$, does
  not exceed $O(A^{s-\theta})$.}
\\[2ex]
Browning and Dietmann \cite{BD} have recently shown that whenever
$s\ge 4$, then
\be{1.9}
\#\{ {\bf a}\in{\cal C}(k,s): \; |{\bf a}| \le A \} \gg A^s,
\ee
so that the estimate in Theorem 1.3 is indeed a non-trivial one. In
particular, it follows that when $s>4k$, then for almost all 
${\bf a} \in{\cal C}(k,s)$,
the equation \rf{1.1} is soluble over $\QQ$.
Since the Hasse principle may fail for ${\bf a}\in{\cal
  C}(k,s)$ only, this implies that the Hasse principle holds for
almost all ${\bf a} \in{\cal C}(k,s)$, but also for almost all ${\bf
  a}\in\ZZ^s$, whenever $s> 4k$. Finally, in the same range for $s$,
Theorem 1.3 implies that for almost all ${\bf a}$ for which
\rf{1.1} has non-trivial integral solutions, there exists a solution
with $0 < |{\bf x}|\le |{\bf a}|^{2/({\hat s}-2k)}$. This last corollary is
rather remarkable, in particular since the upper bound on the size of
the solution is quite small, and not too far from the best possible
one, as the following result shows.
\\[2ex]
{\sc Theorem 1.4.} {\it Let $s > 2k$, and let $\eta >0$. Then, there
  exists a number $c = c(k,s,\eta)>0$ such that the number of ${\bf
    a}$ with $|{\bf a}|\le A$ for which \rf{1.1} admits an integral
  solution in the range $0 < |{\bf x}| \le c|{\bf a}|^{1/(s-k)}$, does
  not exceed $\eta A^s$.}
\medskip

One should compare this with the lower bound \rf{1.9}: even
among the locally soluble equations \rf{1.1}, those that have an
integral solution with $0< |{\bf x}| < c |{\bf a}|^{1/(s-k)}$ form a
thin set, at least when $c$ is small. 
It follows that the exponent $2/({\hat s}-2k)$ that occurs in
Theorem 1.3 cannot be replaced by a number smaller than $1/(s-k)$.

An estimate for the smallest non-trivial solution of an additive
diophantine equation is of considerable importance in diophantine
analysis, also for applications in diophantine approximation; see
Schmidt \cite{Adv} for a prominent example and Birch \cite{70} for
further comments. There are some bounds of this type available in the
literature (eg.\ Pitman \cite{P}), most notably by Schmidt \cite{S1, S2}.
In this context, it is worth recalling that when $s > k^2$ then the
equation \rf{1.1} is soluble over $\QQ_p$, for all primes $p$
(Davenport and Lewis \cite{DL63}). When $k$ is odd, we then expect
that \rf{1.1} is soluble over $\QQ$, and Schmidt \cite{S1} has shown
that for any $\eps >0$ there exists $s_0 (k,\eps)$ such that whenever
$s\ge s_0$ then any equation \rf{1.1} has an integer solution with $0<
|{\bf x}| \ll | {\bf a}|^\eps$. The number $s_0(k,\eps)$ is
effectively computable, but Schmidt's method only yields poor bounds
(see Hwang \cite{Hw} for a discussion of this matter). When $k$ is
even and $s>k^2$, then \rf{1.1} is locally soluble provided only that the $a_j$
are not all of the same sign. In this situation Schmidt \cite{S2}
demonstrated that there still is
some $s_0 (k,\eps)$ such that whenever at least $s_0(k,\eps)$ of the
$a_j$
are positive, and at least $s_0(k,\eps)$ are negative, then
the equation \rf{1.1} is soluble in integers with
$$
0< |{\bf x}| \le |{\bf a}|^{1/k+\eps};
$$
see also  Schlickewei \cite{Schl} when $k=2$.
Schmidt's result is essentially best possible:
if $a\le b$ are coprime natural numbers, and $k$ is even, then any
nontrivial solution of
\be{1.10}
a(x_1^k + \ldots + x_t^k) - b (x_{t+1}^k + \ldots + x_s^k) =0
\ee
must have $b| x_1^k + \ldots + x_t^k$, whence $|{\bf x}| \ge
(b/s)^{1/k}$. Thus, there are equations \rf{1.1} where the smallest
solution is as large as $|{\bf a}|^{1/k}$, even when $s$ is very
large. However, in Theorem 1.3 the exponent $2/({\hat s}-2k)$ is smaller than
$1/k$. It follows that at least when $k$ is even, the exceptional set
for ${\bf a}$ that is estimated in Theorem 1.3, is non-empty. On the
other hand, Theorem 1.3 tell us that examples such as \rf{1.10} 
where the smallest integer
solution is large, must be sparse.
\medskip

We are not aware of any previous attempts to examine additive
diophantine equations on average, save for the recent dissertation of
Breyer \cite{Br}. There, an estimate is obtained that is roughly
equivalent to a variant of Theorem 1.1 in which $B \asymp A^{1/k}$,
and where the sum over ${\bf a}$ is restricted to a rather unnaturally
defined, but reasonably dense subset of $\ZZ^s$. In particular,
Breyer's estimates are not of strength sufficient to derive the Hasse
principle for almost all equations \rf{1.1} with ${\bf a} \in {\cal
  C}(k,s)$, even when $s$ is much larger than $4k$. Yet, our analysis
in section 2
has certain features in common with Breyer's work, most notably the
use of lattice point counts to treat a certain auxiliary equation. The
method
could be described as an attempt to exchange the roles of coefficients and
variables in \rf{1.1}. It is a pure counting device, we cannot
describe the exceptional sets beyond bounds on their cardinality. 
We
postpone a detailed description of our methods until they are needed
in the course of the argument, but remark that the ideas
developed herein can be refined further, and may be applied to related
problems as well. With more work and a different use of the geometry
of numbers, we may advance into the range $3k < s \le 4k$. Perhaps
more importantly, one may derive results similar to those announced as
Theorem 1.3 for the class of general forms of a given degree. Details
must be deferred to sequels of this paper.
\medskip

{\em Notation}. Our notation is standard, or is otherwise explained
within the text. Vectors are typeset in bold, and have dimension $s$
unless indicated otherwise. The symbol ${\bf a}$ is reserved for
tupels
$(a_1,\ldots,a_s)$ with non-zero integers $a_j$. We use
$(x_1;\ldots;x_s)$ to denote the greatest common divisor of the
integers $x_j$. The exponential $\exp(2\pi i \alpha)$ is abbreviated to
$e(\alpha)$. Finally, we apply the familiar $\eps$-convention:
whenever $\eps$ occurs in a statement, it is asserted that the
statement is valid for any positive real number $\eps$. Implicit
constants in
Landau's or Vinogradov's symbols are allowed to depend on $\eps$ in
such circumstances.

\begin{center}
{\bf II. Applications of the geometry of numbers}
\end{center}

\setcounter{section}{2}\setcounter{equation}{0}\noindent
{\bf 2.1. An elementary upper bound estimate.} Our first goal is the
demonstration of Theorem 1.4. The following lattice point count is the main
ingredient.
\\[2ex]
{\sc Lemma 2.1.} {\it Let ${\bf c} \in\ZZ^s$ be a primitive
  vector. Then, for any $X\ge |{\bf c}|$, one has}
$$
\#\{ {\bf x} \in\ZZ^s: \, |{\bf x}|\le X, \, c_1 x_1 + \ldots + c_s
x_s =0 \} \ll |{\bf c}|^{-1} X^{s-1}.
$$
{\it Proof.} See Heath-Brown \cite{Ann}, Lemma 1, for example.
\medskip

Now let $\Xi(A,B)$ denote the number of all ${\bf a}$ with $|{\bf a}|\le
A$ for which the equation \rf{1.1} has an integral solution with $0 <
|{\bf x}| \le B$. We proceed to derive an upper bound for $\Xi(A,B)$.
Note that whenever ${\bf a}$
is counted by $\Xi(A,B)$, then $\rho_{\bf a}(B) -1 \ge 1$. Hence, on
exchanging the order of summation,
$$
\Xi(A,B) \le \sum_{|{\bf a}|\le A} (\rho_{\bf a} (B) -1) = \sum_{0<|{\bf
    x}|\le B} \#\{|{\bf a}| \le A: \; a_1 x_1^k + \ldots
+ a_s x_s^k =0\}
$$
Whenever $B^k\le A$, Lemma 2.1 supplies the estimate
$$
\Xi(A,B) \ll \sum_{0<|{\bf x}|\le B} A^{s-1} \f{(x_1^k ; \ldots ;
  x_s^k)}{|{\bf x}|^k}.
$$
 By symmetry, it suffices
to sum over all ${\bf x}$ with $x_1 = |{\bf x}|$. We sort the
remaining sum according to $d= (x_1; x_2; \ldots; x_s)$. Then $d|x_j$
for all $j$, and we infer that
$$
\Xi(A,B) \ll A^{s-1} \sum_{1\le x_1 \le B} \sum_{d|x_1}
\Big(\f{d}{x_1}\Big)^k \Big(\multsum{y\le x_1}{d|y} 1 \Big)^{s-1}.
$$
Since $x_1/d\ge 1$ holds for all $d|x_1$, it follows that
$$
\Xi(A,B) \ll A^{s-1} \sum_{1\le x\le B} \sum_{d|x}
\Big(\f{x}{d}\Big)^{s-1-k}.
$$
In particular, this confirms the following.
\\[2ex]
{\sc Lemma 2.2.} {\it Let $s\ge k+3$, and suppose that $B^k \le
  A$. Then}
$$
\Xi(A,B) \ll A^{s-1} B^{s-k}.
$$

The proof of Theorem 1.4 is now straightforward. When $s> 2k$ and $0<
C\le 1$, then the choice $B= CA^{1/(s-k)}$ is admissible in Lemma
  2.2. Let $\eta >0$. Then, if $C$ is sufficiently small, Lemma 2.2
  supplies the inequality $\Xi(A,CA^{1/(s-k)}) <\eta A^{s}$. If ${\bf
    a}$ is a vector such that $|{\bf a}|\le A$ and \rf{1.1} has an
  integral solution with $0< |{\bf x}| < C |{\bf a}|^{1/(s-k)}$, then
  ${\bf a}$ is also counted by $\Xi(A, CA^{1/(s-k)})$, and Theorem 1.4
  follows.
\\[2ex]
{\bf 2.2. Another auxiliary mean value estimate.} Our next task is the
derivation of an estimate for the number of solutions of a certain
symmetric diophantine equation. The result will be one of the
cornerstones in the proof of Theorem 1.1. We begin with an
examination of a congruence related to $k$-th powers.
\\[2ex]
{\sc Lemma 2.3.} {\it The number of pairs $(x,y) \in\ZZ^2$ with
  $|x|\le B$, $|y|\le B$ and $x^k \equiv y^k \bmod d$ does not exceed
  $O(B^{1+\eps} + B^{2+\eps} d^{\eps-2/k})$.
\\[2ex]
Proof.} Pairs with $xy=0$ contribute $O(B)$. We sort the remaining
pairs according to the value of $e= (x;y)$, and write $x= e x_0$, $y=
e y_0$. The congruence implies $e^k|d$, and then reduces to $x_0^k
\equiv y_0^k \bmod d e^{-k}$ with $1\le |x_0| \le B/e$, $1\le |y_0|\le
B/e$ and $(x_0; y_0) =1$. There are $2 B/e$ choices for $y_0$, and
since we have now assured that $(x_0; d e^{-k}) =1$, the theory of
$k$-th power residues and a divisor function estimate yield the bound
$O( 1+ B^{1+\eps} d^{-1} e^{k-1})$ for the number of choices for $x_0$,
for any admissible choice of $y_0$. It follows that the number in
question does not exceed
$$
\ll B + B^{2+\eps} \sum_{e^k|d} d^{-1} e^{k-2} + B \sum_{e^k|d}
e^{-1},
$$
which confirms the lemma.
\medskip

Now let $t$ be a natural number, and let $V_t(A,B)$ denote the number
of solutions of the equation
\be{2.1}
\sum_{j=1}^t a_j (x_j^k - y_j^k)=0
\ee
in integers $a_j, x_j, y_j$ constrained to
\be{2.2}
0 < |a_j|\le A, \quad |x_j|\le B, \quad |y_j|\le B, \quad x_j^k \neq
y_j^k.
\ee
{\sc Lemma 2.4.} {\it Let $t\ge 2$, and suppose that $A\ge 2B^k \ge
  1$. Then
$$
V_t(A,B) \ll A^{t-1} (B^{t+1} + B^{2t-k}) B^\eps.
$$
Proof.} We have $|x_j^k - y_j^k| \le 2B^k \le A$. Therefore, by Lemma
2.1,
$$
V_t (A,B) \ll A^{t-1}\underset{1\le j\le t}{\underset{x_j^k \neq
  y_j^k}{  \sum_{|x_j|\le B} \sum_{|y_j|\le B}}}
\f{(x_1^k - y_1^k ; \ldots ; x_t^k - y_t^k)}{\max |x_j^k - y_j^k|}.
$$
By symmetry, it suffices to estimate the portion of the remaining sum
where $|x_j|\le x_1, |y_j|\le x_1$ for all $j$. Then $x_1 >0$, and we
deduce that
$$
V_t(A,B) \ll A^{t-1} \sum_{1\le x_1\le B}
\multsum{|y_1|\le x_1}{y_1^k\neq x_1^k}
\underset{2\le j\le t}{\underset{x_j^k\neq y_j^k}{\sum_{|x_j|\le x_1}
    \sum_{|y_j|\le x_1}}} 
\f{(x_1^k - y_1^k; \ldots ; x_t^k - y_t^k)}{x_1^k - y_1^k}.
$$
For any pair $x_1, y_1$ with $x_1^k \neq y_1^k$, the inner sum will
now be sorted according to the value of $d= (x_1^k - y_1^k ; \ldots;
x_t^k - y_t^k)$. Then $d|x_j^k - y_j^k$ for all $j=1,\ldots,
t$. Therefore, by Lemma 2.3,
\begin{eqnarray*}
V_t (A,B) & \ll &
 A^{t-1} \sum_{1\le x_1\le B} 
\multsum{|y_1|\le  x_1}{y_1^k\neq x_1^k}
\sum_{d|x_1^k - y_1^k}
\f{d}{x_1^k - y_1^k}
\underset{2\le j\le t}{\underset{x_j^k\equiv y_j^k\bmod
    d}{\sum_{|x_j|\le x_1}{ \sum_{|y_j|\le x_1}}}} 1 \\
& \ll &
A^{t-1} \sum_{1\le x_1\le B}
\multsum{|y_1|\le x_1}{y_1^k\neq x_1^k}
\sum_{d|x_1^k - y_1^k}
\f{d}{x_1^k - y_1^k}
(x_1^{1+\eps} + x_1^{2+\eps} d^{-2/k})^{t-1}.
\end{eqnarray*}
Here we apply the trivial inequality $(\xi+\eta)^{t-1} \ll
\xi^{t-1}+\eta^{t-1}$ that is valid for non-negative reals $\xi,\eta$, 
and note that a standard divisor argument
yields
$$
\sum_{1\le x_1 \le B} 
\multsum{|y_1|\le x_1}{y_1^k\neq x_1^k}\sum_{d|x_1^k - y_1^k}
\f{d x_1^{t-1}}{x_1^k - y_1^k}
\ll
B^\eps \sum_{1\le x_1 \le B} x_1^t \ll B^{t+1+\eps}
$$
so that we now deduce that
\be{2.3}
V_t(A,B) \ll A^{t-1}\Big(B^{t+1+\eps} + B^\eps \Upsilon_t(B)\Big)
\ee
with
\be{2.4}
\Upsilon_t(B) = \sum_{1\le x\le B}
\multsum{|y|\le x}{y^k \neq x^k}\sum_{d|x^k - y^k}
\f{x^{2t-2}}{x^k - y^k}
d^{1-2(t-1)/k}.
\ee

We proceed with examining two cases separately. First suppose that
$2(t-1)\ge k$. Then, by a divisor function estimate,
$$
\Upsilon_t(B) \ll B^\eps \sum_{1\le x\le B}
\multsum{|y|\le x}{y^k\neq x^k}
\f{x^{2t-2}}{x^k-y^k}.
$$
When $k$ is even, we group together the two terms $\pm y$, and then
put $h=x-y$. For $y\ge 0$, we have
$$
x^k -y^k = h(x^{k-1} + \ldots + y^{k-1}) \ge hx^{k-1}
$$
whence
\begin{eqnarray*}
\Upsilon_t (B) &\ll&
B^\eps \sum_{1\le x\le B} \sum_{0\le y<x} \f{x^{2t-2}}{x^k -y^k} \\
&\ll&
B^\eps \sum_{1\le x\le B} \sum_{1\le h\le x} h^{-1} x^{2t-1-k} \ll
B^{2t-k+2\eps}.
\end{eqnarray*}
When $k$ is odd, then we first consider the terms with $0\le
y<x$. Then, we may argue as in the case where $k$ is even, and we find
that these pairs $(x,y)$ contribute $O(B^{2t-k+2\eps})$ to
$\Upsilon_t(B)$. The remaining terms, with $-x\le y<0$, are even simpler
to control. Since $k$ is odd, we have $x^k - y^k\ge x^k$, and so,
$$
\sum_{1\le x\le B} \sum_{-x\le y<0} \f{x^{2t-2}}{x^k-y^k} \le
\sum_{1\le x\le B} x^{2t-1-k} \ll B^{2t-k+\eps}.
$$
It follows that $\Upsilon_t(B) \ll B^{2t-k+\eps}$ holds in all cases, and
by \rf{2.3}, we have now shown that whenever $2(t-1)\ge k$, one
has
$$
V_t(A,B) \ll A^{t-1} (B^{t+1+\eps} + B^{2t-k+\eps}),
$$
as required.
\medskip

It remains to investigate the situation where $2(t-1)< k$. Here,
a divisor function estimate applied within \rf{2.4} yields
$$
\Upsilon_t(B) \ll \sum_{1\le x\le B} x^{2t-2}
\multsum{|y|\le x}{y^k\neq x^k}
(x^k-y^k)^{-2(t-1)/k}.
$$
When $k$ is even, we manipulate this sum much as in the previous case,
and find that
\begin{eqnarray*}
\Upsilon_t(B) &\ll &
\sum_{1\le x\le B} x^{2t-2}\sum_{0\le y <x}(x^k -
y^k)^{-2(t-1)/k}\\
&\ll &
\sum_{1\le x\le B} x^{2t-2 -2(k-1)(t-1)/k } \sum_{1\le h\le x}
h^{-2(t-1)/k} 
\ll 
B^{2+\eps}.
\end{eqnarray*}
A similar computation yields the same result when $k$ is
odd. Therefore, when $2(t-1)<k$, we now deduce from \rf{2.3}
that
$$
V_t (A,B) \ll A^{t-1} (B^{t+1+\eps} + B^{2+\eps})\ll A^{t-1}
B^{t+1+\eps}.
$$
This confirms the claim in Lemma 2.4.
\medskip

Now let $U_t (A,B)$ denote the number of solutions of \rf{2.1} in
integers $a_j, x_j, y_j$ satisfying
$$
0<|a_j|\le A, \quad |x_j|\le B, \quad |y_j|\le B.
$$
For any solution counted by $U_t(A,B)$, let $r$ be the number of
$j\in\{1,2,\ldots, t\}$ with $x_j^k\neq y_j^k$. The contribution to
$U_t (A,B)$ made by solutions with $r=0$ is obviously no larger than
$O(A^t B^t)$. By symmetry, we now deduce that
$$
U_t (A,B) \ll A^t B^t + \sum_{r=1}^t (AB)^{t-r} V_r (A,B).
$$
The definition of $V_1(A,B)$ implies that $V_1(A,B) =0$. We now
suppose that $A\ge 2B^k \ge 1$, apply Lemma 2.4 to bound $V_r(A,B)$
for $2\le r\le t$, and then deduce the 
 following estimate.
\\[2ex]
{\sc Theorem 2.5.} {\it Let $t\ge 2$. Then, for real numbers $A,B$
  with $A\ge 2B^k \ge 1$, one has}
$$
U_t (A,B) \ll (AB)^t + A^{t-1} B^{2t-k+\eps}.
$$
\bigskip

\begin{center}
{\bf III. Local solubility}
\end{center}

\setcounter{section}{3}\setcounter{equation}{0}\noindent
{\bf 3.1. The singular integral.} Local solubility of additive
equations has been investigated by Davenport and Lewis \cite{DL63},
and by Davenport \cite{notes}. The analytic condition \rf{1.5} for
local solubility is implicit in \cite{DL63}. Unfortunately, these
prominent references are insufficient for our purposes. A lower bound
for $J_{\bf a}{\fk S}_{\bf a}$ in terms of $|{\bf a}|$ is needed
whenever this product in non-zero, at least for almost all ${\bf
  a}$. An estimate of this type is supplied in this section.
\medskip

We begin with the singular integral. Most of our work is routine, so
we shall be brief. When $\beta\in\RR$, $B>0$, let 
\be{3.1}
v(\beta,B)= \int_{-B}^B e(\beta \xi^k) d\xi.
\ee
A partial integration readily confirms the bound
\be{3.2}
v(\beta, B) \ll B(1+B^k |\beta|)^{-1/k}
\ee
whence whenever $s>k$ one has
\be{3.3}
\int_{-\inf}^\inf |v(\beta,B)|^s d\beta \ll B^{s-k}.
\ee
We also see that for $s>k$ and ${\bf a}\in(\ZZ\bs \{0\})^s$, the
integral
\be{3.4}
J_{\bf a}(B) =\int_{-\inf}^\inf v(a_1\beta , B)\ldots v(a_s \beta,B)
d\beta
\ee
converges absolutely. By H\"older's inequality and \rf{3.3},
\begin{eqnarray}\lefteqn{
\int_{-\inf}^\inf |v(a_1\beta, B) \ldots v(a_s \beta, B)| d\beta}
\nonumber \\
&\le &
\prod_{j=1}^s \Big( \int_{-\inf}^\inf |v(a_j\beta, B)|^s d\beta
\Big)^{1/s}
\ll
|a_1 \ldots a_s|^{-1/s} B^{s-k}. \label{3.5}
\end{eqnarray}
In particular, it follows that
\be{3.6}
J_{\bf a} (B) \ll |a_1 \ldots a_s|^{-1/s} B^{s-k}.
\ee

The integral $J_{\bf a}(B)$ arises naturally as the singular integral
in our application of the circle method in section 4. The dependence
on $B$ can be made more explicit. By \rf{3.1}, one has $v(\beta, B)=
Bv (\beta B^k, 1)$. Now substitute $\beta$ for $\beta B^k$ in \rf{3.4}
to infer that
\be{3.7}
J_{\bf a}(B) = B^{s-k} J_{\bf a}
\ee
where $J_{\bf a} = J_{\bf a}(1)$ is the number that occurs in
\rf{1.2}, and in Theorem 1.1.
\medskip

It remains to establish a lower bound for $J_{\bf a}$. The argument
depends on the parity of $k$, and we shall begin with the case when
$k$ is even. Throughout, we suppose that
\be{3.8}
|a_s|\ge |a_j| \quad (1\le j<s).
\ee
Define $\sigma_j= a_j/|a_j|\in \{ 1,-1\}$. Then, by \rf{3.1},
$$
v(a_j \beta,1) = 2\int_0^1 e(a_j\beta \xi^k) d\xi =
\f{2}{k} |a_j|^{-1/k} \int_0^{|a_j|} \eta^{(1-k)/k} e(\sigma_j \beta
\eta) d\eta.
$$
Let ${\fk A} = [0,|a_1|] \times \ldots \times [0,|a_s|]$, and define
the linear form $\tau$ through the equation
\be{3.9}
\sigma_s \tau = \sigma_1 \eta_1 + \ldots + \sigma_s \eta_s.
\ee
Then, we may rewrite \rf{3.4} as
$$
J_{\bf a} = \Big(\f{2}{k}\Big)^s |a_1 \ldots a_s|^{-1/k}
\int_{-\inf}^\inf \int_{\fk A} (\eta_1 \ldots \eta_s)^{(1-k)/k}
e(\sigma_s \tau \beta) d\ETA d\beta.
$$
Now substitute $\tau$ for $\eta_s$ in the innermost integral. Then, by
Fubini's theorem and \rf{3.9},
$$
\int_{\fk A} (\eta_1 \ldots \eta_s)^{(1-k)/k} e(\sigma_s \tau \beta)
d\ETA = \int_{-\inf}^\inf E(\tau) e(\sigma_s \tau\beta) d\tau
$$
where
\be{3.9a}
E(\tau) = \int_{{\fk E}(\tau)} (\eta_1 \ldots \eta_{s-1} \eta_s(\tau,
\eta_1, \ldots, \eta_{s-1}))^{(1-k)/k} d(\eta_1, \ldots, \eta_{s-1}),
\ee
in which $\eta_s$ is the linear form defined implicitly by \rf{3.9},
and where ${\fk E}(\tau)$ is the set of all $(\eta_1, \ldots,
\eta_{s-1})$ satisfying the inequalities
\begin{eqnarray*}
0 &\le& \eta_j \le |a_j| \qquad (1\le j<s), \\
0 &\le & \tau- \sigma_s \sigma_1 \eta_1 - \sigma_s \sigma_2 \eta_2
-\ldots - \sigma_s \sigma_{s-1} \eta_{s-1}\le |a_s|.
\end{eqnarray*}
It transpires that $E$ is a non-negative continuous function with
compact support, and that for $\tau$ near $0$, this function is of
bounded variation. Therefore, by Fourier's integral theorem,
$$
\lim_{N\to\inf} \int_{-N}^N \int_{-\inf}^\inf E(\tau) e(\sigma_s
\tau\beta) d\tau d\beta = E(0),
$$
and we infer that
\be{3.10}
J_{\bf a} = \Big(\f{2}{k}\Big)^s |a_1 \ldots a_s|^{-1/k} E(0).
\ee
In particular, it follows that $J_{\bf a} \ge 0$. Also, when all $a_j$
have the same sign, then ${\fk E}(0)= \{{\bf 0}\}$, and \rf{3.10}
yields $J_{\bf a}=0$.
\medskip

Now suppose that not all the $a_j$ are of the same sign. First,
consider the situation where $\sigma_1 = \ldots = \sigma_{s-1}$. Then
we have $\sigma_s \sigma_j = -1$ $(1\le j<s)$. By \rf{3.8}, we see
that the set of $(\eta_1, \ldots, \eta_{s-1})$ defined by
$$
\f{|a_j|}{2s} \le \eta_j \le \f{|a_j|}{s}
\qquad (1\le j<s)
$$
is contained in ${\fk E}(0)$, and its measure is bounded
below by $(2s)^{-s}|a_1 a_2 \ldots a_{s-1}|$. By \rf{3.9a}, we now
deduce that
$$
E(0)\gg |a_1 \ldots a_{s-1}|^{1/k} |a_s|^{(1-k)/k},
$$
and \rf{3.10} then implies the bound $J_{\bf a} \gg |a_s|^{-1} = |{\bf
  a}|^{-1}$.
\medskip

In the remaining cases, both signs occur among $\sigma_1, \ldots,
\sigma_{s-1}$. We may therefore suppose that for some $r$ with $2\le r
< s$ we have
$$
\sigma_s \sigma_j = -1 \quad (1\le j<r), \qquad \sigma_s \sigma_j =1
\quad (r\le j<s).
$$
Take $\tau=0$ in \rf{3.9}. Then $\eta_s$ is the linear form
\be{3.11}
\eta_s = \eta_1 + \ldots + \eta_{r-1} - \eta_r - \ldots - \eta_{s-1}.
\ee
By symmetry, we may suppose that
$$
|a_1|\le |a_2|\le \ldots \le |a_{r-1}|, \qquad |a_r|\le |a_{r+1}| \le
\ldots \le |a_{s-1}|.
$$
We define $t$ by $t= r-1$ when $|a_{r-1}|\le |a_r|$, and otherwise as
the largest $t$ among $r, r+1, \ldots, s-1$ where $|a_t|\le
|a_{r-1}|$. Now consider the set of $(\eta_1, \ldots, \eta_{s-1})$
defined by the inequalities
\begin{eqnarray*}
\f{|a_j|}{2s} &\le  \eta_j \le & \f{|a_j|}{s} \qquad (1\le j\le r-1),
\\
\f{|a_j|}{8s^2} &\le \eta_j\le & \f{|a_j|}{4s^2} \qquad (r\le j\le
t),\\
\f{|a_{r-1}|}{8s^2} &\le \eta_j \le & \f{|a_{r-1}|}{4s^2} \quad 
(t<j<s).
\end{eqnarray*}
It is readily checked that on this set, the number $\eta_s$ defined in
\rf{3.11} satisfies the inequalities $\f{|a_{r-1}|}{4s} \le \eta_s \le
|a_{r-1}|$. Moreover, the measure of this set is $\gg |a_1 \ldots a_t|
|a_{r-1}|^{s-t+2}$. By \rf{3.9a}, it follows that
$$
E(0) \gg |a_1 \ldots a_t|^{1/k} |a_{r-1}|^{(s-t+2)/k}
|a_{r-1}|^{(1-k)/k},
$$
and again one then deduces from \rf{3.10} the bound $J_{\bf a} \gg
|{\bf a}|^{-1}$.
\medskip

Finally, we discuss the case where $k$ is odd. The main differences in
the treatment occur in the initial steps. When $k$ is odd, one may
transform \rf{3.1} into
$$
v(a_j \beta, 1) = \f{1}{k} |a_j|^{-1/k} \int_0^{|a_j|} \eta^{(1-k)/k}
(e(\beta\eta) + e(-\beta \eta)) d\eta.
$$
Let ${\SIGMA}= (\sigma_1, \ldots, \sigma_s)$ with $\sigma_j \in
\{1,-1\}$. For any such ${\SIGMA}$, define $\tau$ through
\rf{3.9}. Then, following through the argument used in the even case,
we first arrive at the identity
$$
J_{\bf a} = k^{-s} |a_1 \ldots a_s|^{-1/k} \sum_{{\boldmath\mbox{$\scriptstyle\sigma$}}}
\int_{-\inf}^\inf \int_{\fk A} (\eta_1\ldots \eta_s)^{(1-k)/k}
e(\sigma_s \tau \beta) d\ETA d\beta.
$$
Here the sum is over all $2^s$ choices of ${\SIGMA}$. Again as
before, we see that each individual summand is non-negative, and when
not all of $\sigma_1,\ldots, \sigma_s$ have the same sign, then one
finds the lower bound $\gg |{\bf a}|^{-1}$ for this summand. Thus, we
now see that $J_{\bf a}\gg |{\bf a}|^{-1}$ again holds, this time for
any choice of ${\bf a}$.
\medskip\noindent

For easy reference, we summarize the above results as a lemma.
\\[2ex]
{\sc Lemma 3.1.} {\it Suppose that $s>k$. Then the singular integral
  $J_{\bf a}$ converges absolutely, and one has $0\le J_{\bf a}\ll
  |a_1 a_2 \ldots a_s|^{1/s}$. Furthermore, when $k$ is odd, or when $k$ is
  even and $a_1, \ldots, a_s$ are not all of the same sign, then
 $J_{\bf
    a}\gg |{\bf a}|^{-1}$. Otherwise
  $J_{\bf a}=0$.}
\\[2ex]
{\bf 3.2. The singular series.} In the introduction, we defined the
classical singular series as a product of local densities. We briefly
recall its representation as a series. Though this is standard in
principle, our exposition makes the dependence on the coefficients
${\bf a}$ in \rf{1.1} as explicit as is necessary for the proof of
Theorem 1.2 in the next section. Recall that $k\ge 3$.
\smallskip

For $q\in\NN$, $r\in\ZZ$ define the Gaussian sum
\be{3.20}
S(q,r) = \sum_{x=1}^q e(rx^k/q).
\ee
Let $\kappa(q)$ be the multiplicative function that, on prime powers
$q=p^l$, is given by
$$
\kappa (p^{uk+v}) = p^{-u-1} \quad (u\ge 0, \, 2\le v \le k), \quad
\kappa(p^{uk+1}) = k p^{-u-{1}/{2}}.
$$
Then, as a corollary to Lemmas 4.3 and 4.4 of Vaughan \cite{hlm}, one
has $S(q,r) \ll q\kappa (q)$ whenever $(q;r) =1$, and one concludes
that
\be{3.21}
q^{-1} S(q,r) \ll \kappa (q/(q;r))
\ee
holds for all $q\in\NN, r\in\ZZ$. Now let
\be{3.22}
T_{\bf a}(q) = q^{-s} \multsum{r=1}{(r;q)=1}^q S(q,a_1r)\ldots
S(q,a_sr).
\ee
Then, by \rf{3.21}, 
\be{3.23}
T_{\bf a}(q) \ll q\kappa (q/(q; a_1))\ldots \kappa (q/(q;a_s)).
\ee
Moreover, by working along the proof of Lemma 2.11 of Vaughan
\cite{hlm}, one finds that $T_{\bf a}(q)$ is a multiplicative function
of $q$. Also, one can use the definition of $\kappa$ to confirm
that whenever $s\ge k+2$ then the expression on the right hand side of
\rf{3.23} may be summed over $q$ to an absolutely convergent
series. Thus, we may also sum $T_{\bf a}(q)$ over $q$ and rewrite the
series as an Euler product. This gives
\be{3.24}
\sum_{q=1}^\inf T_{\bf a} (q) = \prod_p \sum_{h=0}^\inf T_{\bf a}
(p^h).
\ee
However, by \rf{3.20} and \rf{3.22}, and orthogonality,
\be{3.25}
\sum_{h=0}^l T_{\bf a} (p^h) = p^{-ls} \sum_{r=1}^{p^l} S(p^l, r a_1)
\ldots S(p^l, r a_s) = p^{l(1-s)} M_{\bf a}(p^l)
\ee
where $M_{\bf a}(p^l)$ is the number of incongruent solutions of the
congruence
$$
a_1 x_1^k + \ldots + a_s x_s^k \equiv 0 \bmod p^l.
$$
We may take the limit for $l\to\inf$ in \rf{3.25} because all sums in
\rf{3.24} are convergent. This shows that the limit $\chi_p$, as
defined in \rf{1.3}, exists. In view of \rf{3.24} and \rf{1.4}, we may
summarize our results as follows.
\\[2ex]
{\sc Lemma 3.2.} {\it Let $s\ge k+2$. Then, for any ${\bf
    a}\in(\ZZ \bs \{0\})^s$, the singular product \rf{1.4} converges,
  and has the alternative representation}
$$
{\fk S}_{\bf a} = \sum_{q=1}^\inf T_{\bf a}(q).
$$

A slight variant of the preceding argument also supplies an
estimate for $\chi_p ({\bf a})$ when $p$ is large.
\\[2ex]
{\sc Lemma 3.3.} {\it Let $s\ge k+2$. Then there is a real number
  $c=c(k,s)$ such that for any choice of $a_1, \ldots, a_s \in\ZZ\bs\{
  0\}$ for which at least $k+2$ of the $a_j$ are not divisible by $p$,
  one has}
$
|\chi_p ({\bf a}) -1| \le cp^{-2}.
$
\\[1ex]
{\it Proof.} We begin with \rf{3.25}, and note that $T_{\bf a} (1) =
1$. Then
$$
p^{l(1-s)} M_{\bf a} (p^l) -1 = \sum_{h=1}^l T_{\bf a} (p^h).
$$
One has $\kappa(q) \le k$ for any prime power $q$. Hence, by
\rf{3.23}, and since $k+2$ of the $a_j$ are coprime to $p$, one
finds that $|T_{\bf a}(p^h) |\le k^s \kappa (p^h)^{k+2}
p^h$. Consequently, a short calculation based on the definition of
$\kappa$ reveals that
$$
|p^{l(1-s)} M_{\bf a}(p^l) -1| \le
k^s \sum_{h=1}^l \kappa (p^h)^{k+2} p^h
\le 
k^{s+k+2} p^{-2}.
$$
The lemma follows on considering the limit $l\to\inf$.
\\[2ex]
{\bf 3.3. Proof of Theorem 1.2.} Throughout, we suppose
that $s\ge k+3$.
 For ${\bf a}\in(\ZZ\bs\{0\})^s$, let ${\cal
  S}({\bf a})$ denote the set of all primes that divide at least two
of the integers $a_j$. Lemma 3.3 may then be applied to all primes
$p\notin {\cal S}({\bf a})$, and we deduce that there exists a number
$C=C(k,s)>0$ such that the inequalities
\be{3.30}
\f{1}{2} \le \multprod{p\notin {\cal S}({\bf a})}{p>C}
\chi_p({\bf a}) \le 2
\ee
hold for all ${\bf a}$.
It will be convenient to write
$$
{\cal P}({\bf a}) = {\cal S}({\bf a}) \cup \{p: p\le C\};
$$
this set contains all primes not covered by \rf{3.30}. For a prime
$p\in{\cal P}({\bf a})$, let
$$
l(p) = \max \{l: p^l|a_j \mbox{ for some } j\},
$$
and then define the numbers
$$
P({\bf a}) = \prod_{p\in{\cal P}({\bf a})} p, \quad
P_0({\bf a}) = \multprod{p\in{\cal S}({\bf a})}{p>C} p, \quad
P^{\dagger} ({\bf a}) = \prod_{p\in{\cal P}({\bf a})} p^{l(p)}.
$$
For later use, we note that
$$
P_0({\bf a}) | P({\bf a}), \quad P({\bf a}) | H P_0({\bf a})
$$
in which we wrote
$$
H= \prod_{p\le C} p.
$$

Now fix a number $\del >0$, to be determined later, and consider the
sets
\begin{eqnarray}
{\cal A}_1&  = &\{ |{\bf a}|\le A: \, P({\bf a}) > A^\del\}, \label{3.31} \\
{\cal A}_2& = & \{|{\bf a}|\le A: \, P({\bf a}) \le A^\del, \,
P^\dagger ({\bf a}) > A^{2\del}\} \label{3.32}
\end{eqnarray}
It transpires that the set ${\cal A}_1 \cup {\cal A}_2$ contains all
${\bf a}$ where the singular series is likely to be
smallish. Fortunately, ${\cal A}_1$ and ${\cal A}_2$ are defined by
divisibility constraints that are related to convergent sieves, so one
expects ${\cal A}_1, {\cal A}_2$ to be thin sets. This is indeed the
case, as we shall now show.

We begin by counting elements of ${\cal A}_1$. For a natural number
$d$, let ${\cal A}_1(d)= \{{\bf a}\in {\cal A}_1: P_0({\bf a})
=d\}$. If there is some ${\bf a}
\in{\cal A}_1(d)$, then by the definition of ${\cal S}({\bf a})$, we have
$d^2| a_1 a_2 \ldots a_s$, whence $d\le A^{s/2}$. On the other hand,
$A^\del < P({\bf a}) \le H P_0 ({\bf a}) \le H d$. This shows that
$$
\#{\cal A}_1 =  \sum_{A^\del / H < d \le A^{s/2}} \# {\cal A}_1 (d) 
\le\sum_{A^\del / H < d \le A^{s/2}} \#\{|{\bf a}|\le A: d^2| a_1
\ldots a_s\}.
$$
By a standard divisor argument, we may conclude that
\be{3.33}
\#{\cal A}_1 \ll A^{s+\eps} \sum_{A^\del/H < d\le A^{s/2}} d^{-2} \ll
A^{s-\del+\eps}.
\ee

The estimation of $\#{\cal A}_2$ proceeds along the same lines, but we
will have to bound the number of integers with small square-free
kernel. When $n$ is a natural number, let
$$
n^* = \prod_{p|n} p
$$
denote its squarefree kernel. One then has the following simple bound
(Tenenbaum \cite{T}, Theorem II.1.12).
\\[2ex]
{\sc Lemma 3.4.} {\it Let $\nu\ge 1$ be a real number. Then,}
$$
\# \{ n\le X^\nu: n^* \le X\} \ll X^{1+\eps}.
$$
For $d\in\NN$, let ${\cal A}_2(d) = \{{\bf a}\in{\cal A}_2 : P^\dagger
({\bf a}) = d \}$. Since we have $P^\dagger({\bf a}) | a_1 a_2 \ldots a_s$,
we must have
$$
A^{2\del} < d\le A^s
$$
whenever ${\cal A}_2(d)$ is non-empty. Moreover, $P({\bf a})$ is the
square-free kernel of $P^\dagger ({\bf a})$, so that $d^* \le
P^\del$. This yields the bound
$$
\#{\cal A}_2 =
\multsum{A^{2\del} < d\le A^s}{d^* \le A^\del} \# {\cal A}_2 (d) 
\le 
\multsum{A^{2\del} < d\le A^s}{d^* \le A^\del} \# \{|{\bf a}| \le A:
d| a_1 a_2 \ldots a_s\}.
$$
The divisor argument used within the estimation of $\#{\cal A}_1$ also
applies here, and gives
$$
\#{\cal A}_2 \ll A^{s+\eps} \multsum{A^{2\del}<d \le A^s}{d^*\le
  A^\del}
\f{1}{d} \ll A^{s-2\del+\eps}
\multsum{d\le A^s}{d^*\le A^\del} 1.
$$
By Lemma 3.4, it follows that
\be{3.34}
\#{\cal A}_2 \ll A^{s-\del + \eps}.
\ee

We are ready to establish Theorem 1.2. It will suffice to find a lower
bound for ${\fk S}_{\bf a}$ for those $|{\bf a}| \le A$ where ${\fk
  S}_{\bf a}>0$ and ${\bf a}\notin {\cal A}_1\cup{\cal A}_2$. Let
$p\in{\cal P}({\bf a})$. We have $\chi_p({\bf a})>0$, whence
\rf{1.1} is soluble in $\QQ_p$. By homogeneity, there is then a
solution ${\bf x}\in\ZZ_p$ of \rf{1.1} with $p\nmid {\bf x}$. In
particular, for any $h\in\NN$, we can find integers $y_1, \ldots, y_s$
that are not all divisible by $p$, and satisfy the congruence
\be{3.40}
a_1 y_1^k + \ldots + a_s y_s^k \equiv 0 \bmod p^h.
\ee
It will be convenient to rearrange indices to assure
that $p\nmid y_1$. Let $\nu(p)$ be defined by $p^{\nu(p)}\| k$, and
recall that a $k$-th power residue $\bmod p^{\nu(p)+2}$ is also a
$k$-th power residue modulo $p^\nu$, for any $\nu\ge \nu(p)+2$. We
choose $h= l(p) + \nu(p)+2$ in \rf{3.40}, and define $e$ by
$p^e\|a_1$. For $l>h$, choose numbers $x_j$, for $2\le j\le s$, with
$1\le x_j \le p^l$ and $x_j\equiv y_j \bmod p^h$. Then, by \rf{3.40},
$$
-\f{a_1}{p^e} y_1^k \equiv \f{a_2 x_2^k + \ldots + a_s x_s^k}{p^e}
\bmod p^{h-e},
$$
and we have $e\le l(p)$, whence $h-e \ge \nu(p) +2$. Thus, for any
choice of $x_2, \ldots, x_s$ as above, there is a number $x_1$ with
$$
a_1 x_1^k + \ldots + a_s x_s^k \equiv 0 \bmod p^l.
$$
Counting the number of possibilities for $x_2, \ldots, x_s$ yields
$M_{\bf a} (p^l) \ge p^{(s-1)(l-h)}$, and consequently,
$$
\chi_p ({\bf a}) \ge p^{(1-s)h}.
$$
We may combine this with \rf{3.30} to infer that
\be{3.410}
{\fk S}_{\bf a} \ge \f{1}{2} \prod_{p\in{\cal P}({\bf a})} p^{(1-s)h}.
\ee
In this product, we first consider primes $p\in{\cal P}({\bf a})$ where
$l(p)=0$. Then $p\nmid a_1 a_2 \ldots a_s$, and the definition of
${\cal P}({\bf a})$ implies that $p\le C$. Also, since $\nu(p)\le k$,
we have $h\le k+2$ so that
$$
\multprod{p\in{\cal P}({\bf a})}{l(p)=0} p^{(1-s)h} \ge
\prod_{p\le C} p^{(1-s)(k+2)} \ge H^{(1-s)(k+2)}.
$$
Next, consider $p\in{\cal P}({\bf a})$ with $l(p)\ge 1$. Then, much as
before, $h\le k+2 + l(p) \le l(p) (k+3)$. Hence,
$$
\multprod{p\in{\cal P}({\bf a})}{l(p)\ge 1} 
p^{(1-s)h} \ge P^\dagger({\bf a})^{(1-s)(k+3)}.
$$
However, since ${\bf a}\notin {\cal A}_1\cup {\cal A}_2$, we have
$P^\dagger ({\bf a}) \le A^{2\del}$, so that we now deduce from
\rf{3.410} that
\be{3.42}
{\fk S}_{\bf a} \gg A^{2\del(1-s)(k+3)}.
\ee
The synthesis is straightforward. Let $\gm >0$. Then choose $\del
=\gm /(8(s-1)(k+3))$, and suppose that $A$ is large. Then
\rf{3.42} implies that ${\fk S}_{\bf a} > A^{-\gm}$. If that fails,
then ${\fk S}_{\bf a}=0$, or else ${\bf a}\in{\cal A}_1 \cup {\cal
  A}_2$. The estimates \rf{3.33} and \rf{3.34} imply Theorem 1.2.
\\[2ex]
{\bf 3.4. An auxiliary upper bound.} We close this section with a
succession of lemmata that involve the function $\kappa$, 
and that will provide an upper bound for ${\fk S}_{\bf a}$
on average. The results will be relevant for the application of the
circle method in the next section.
\\[2ex]
{\sc Lemma 3.5.} {\it One has}
$$
\sum_{d|q} d\kappa (d) \ll q^{1+\eps} \kappa(q).
$$
{\it Proof.} Let $p$ be a prime, and suppose that $0\le j\le l$. Then,
an inspection of the definition of $\kappa$ readily reveals that the
crude inequality $p^j \kappa (p^j) \le k p^l \kappa (p^l)$
holds. Consequently, one also has
$$
\sum_{d|p^l} d\kappa (d) = \sum_{j=1}^l p^j \kappa (p^j) \le k(l+1)
\kappa (p^l) p^l.
$$
By multiplicativity, this implies the bound
$$
\sum_{d|q} d\kappa (d) \le q\kappa (q) \prod_{p^l\|q} k(l+1),
$$
which is more than required.
\\[2ex]
{\sc Lemma 3.6.} {\it Uniformly for $q\in\NN$ and $A\ge 1$, one has}
$$
\sum_{1\le a\le A} \kappa(q/(q;a)) \ll A q^\eps \kappa (q).
$$
{\it Proof.} We sort the $a$ according to the value of $d=
(q;a)$. Then
\begin{eqnarray*}
\sum_{1\le a\le A} \kappa (q/(q;a)) &=&
\multsum{d|q}{d\le A} \kappa (q/d) \multsum{1\le a\le A}{(a;q)=d} 1 \\
&\le &
A \sum_{d|q} d^{-1} \kappa (q/d) = \f{A}{q} \sum_{d|q} d\kappa (d).
\end{eqnarray*}
The lemma now follows by appeal to Lemma 3.5.
\\[2ex]
{\sc Lemma 3.7.} {\it Let $s\ge k+2$. Then}
$$
\sum_{|{\bf a}|\le A} \sum_{q=1}^\inf (a_1 a_2 \ldots a_s)^{-1/s} q^{
  1+1/(2k)} \kappa (q/(q;a_1) ) \ldots \kappa (q/(q;a_s)) \ll A^{s-1}.
$$
{\it Proof.} The terms to be summed are non-negative. Thus, we may
take the sum over ${\bf a}$ first. This then factorizes, and by Lemma
3.6 and partial summation, the left hand side in Lemma 3.7 is seen not
to exceed
$$
A^{s-1} \sum_{q=1}^\inf q^{1+1/(2k)+\eps} \kappa (q)^s.
$$
The remaining sum converges for $s\ge k+2$, as one readily confirms by
considering the corresponding Euler product. The lemma follows.
\medskip

We now apply the last estimate to the singular series. Let $T_{\bf
  a}(q)$ be as in \rf{3.22}. When $Q\ge 1$, define the tail of ${\fk
  S}_{\bf a}$ as
\be{3.41}
{\fk S}_{\bf a}(Q) = \sum_{q\ge Q} T_{\bf a}(q)
\ee
which is certainly convergent for $s\ge k+2$; compare Lemma 3.2. Also,
note that ${\fk S}_{\bf a}= {\fk S}_{\bf a}(1)$.
\\[2ex]
{\sc Lemma 3.8.} {\it Let $s\ge k+2$. Then, uniformly in $A\ge 1$,
  $Q\ge 1$, one has}
$$
\sum_{|{\bf a}| \le A} (a_1 a_2 \ldots a_s)^{-1/s} |{\fk S}_{\bf
  a}(Q)|\ll A^{s-1} Q^{-1/(2k)}.
$$
{\it Proof.} By \rf{3.23},
\begin{eqnarray*}
|{\fk S}_{\bf a}(Q)| &\le &
Q^{-1/(2k)} \sum_{q=1}^\inf q^{1/(2k)} |T_{\bf a}(q)|\\
&\le &
Q^{-1/(2k)} \sum_{q=1}^\inf q^{1+1/(2k)} \kappa (q/(q;a_1)) \ldots
  \kappa(q/(q;a_s)),
\end{eqnarray*}
and the lemma follows from Lemma 3.7.
\bigskip

\begin{center}
{\bf IV. The circle method}
\end{center}

\setcounter{section}{4}\setcounter{equation}{0}
{\bf 4.1. Preparatory steps.} In this section, we establish Theorem
1.1.
The argument is largely standard, save for the ingredients to be
imported from the previous sections of this memoir.
\medskip

We employ the following notational convention throughout this section:
if $h: \RR \to \CC$ is a function, and ${\bf a}\in\ZZ^s$, then we
define
\be{4.1}
h_{\bf a} (\al) = h(a_1\al) h(a_2\al) \ldots h(a_s\al).
\ee
As is common in problems of an additive nature, the Weyl sum
\be{4.2}
f(\al) = \sum_{|x|\le B} e(\al x^k)
\ee
is prominently featured in the argument to follow, because by
orthogonality, one has
\be{4.3}
\rho_{\bf a} (B) = \int_0^1 f_{\bf a} (\al) d\al.
\ee
The circle method will be applied to the integral in \rf{4.3}. With
applications in mind that go well beyond those in the current
communication, we shall treat the ``major arcs'' under very mild
conditions on $A,B$, and for the range $s\ge k+2$.
\medskip

Let $A\ge 1$, $B\ge 1$, and fix a real number $\eta >0$. Then put
$Q=B^\eta$. Let ${\fk M}$ denote the union of the intervals
\be{4.4}
\Big|\al - \f{r}{q}\Big| \le \f{Q}{AB^k}
\ee
with $1\le r\le q <Q$, and $(r,q)=1$. When $\eta \le \f{1}{3}$, these
intervals are pairwise disjoint, and we write ${\fk m}=[Q/(AB^k), 1+Q
/(AB^k)]\bs {\fk M}$. When ${\fk A}$ is one of ${\fk M}$ or ${\fk m}$,
let
\be{4.5}
\rho_{\bf a} (B,{\fk A}) = \int_{\fk A}f_{\bf a} (\al) d\al
\ee
and note that
\be{4.6}
\rho_{\bf a} (B) = \rho_{\bf a} (B,{\fk M}) + \rho_{\bf a}(B,{\fk m}).
\ee
\\[2ex]
{\bf 4.2. The major arc analysis.} In this section we make heavy use
of the results in Vaughan's book \cite{hlm} on the subject. He works
with the Weyl sum
$$
g(\al) = \sum_{1\le x\le B} e(\al x^k)
$$
that is related with our $f$ through the formulae
$$
f(\al) = 1+2g (\al) \quad (k \mbox{ even}), \quad f(\al) = 1+g(\al) +
g(-\al) \quad (k \mbox{ odd}).
$$
Thus, in particular, Theorem 4.1 of \cite{hlm} yields the following.
\\[2ex]
{\sc Lemma 4.1.} {\it Let $\al \in \RR$, $r\in\ZZ$, $q\in\NN$ and
  $a\in\ZZ$ with $a\neq 0$. Then}
$$
f(a\al) = q^{-1} S(q,ar) v(a(\al - r/q)) + O(q^{1/2 +\eps} 
(1 + |a| B^k |\al - r/q|)^{1/2}).
$$
Here, and throughout the rest of this section, we define 
$v(\beta) = v(\beta, B)$ through \rf{3.1}. 
When $|a|\le A$, and $\al \in{\fk M}$ is in the interval \rf{4.4}, we
find that
$$
f(a\al) = q^{-1} S(q,ar) v(a(\al-r/q)) + O(Q^2).
$$
This we use with $a=a_j$ and multiply together. Then
$$
f_{\bf a} (\al) = q^{-s} S(q,a_1 r) \ldots S(q,a_s r) v_{\bf a} (\al
-r/q) + O(Q^2 B^{s-1}).
$$
Now integrate over ${\fk M}$, and recall the definition of the
latter. By \rf{4.5} and \rf{3.22}, we then arrive at
$$
\rho_{\bf a}(B,{\fk M}) = \sum_{q< Q} T_{\bf a}(q)
\int_{-Q/(AB^k)}^{Q/(AB^k)} v_{\bf a} (\beta) d\beta + O(Q^5 B^{s-1-k}
A^{-1}).
$$
Here, we complete the sum over $q$ to the singular series, and the
integral over $\beta$ to the singular integral. On writing
\be{4.7}
\int_{-Q/(AB^k)}^{Q/(AB^k)} v_{\bf a} (\beta) d\beta = J_{\bf a} (B) +
E_{\bf a},
\ee
we may recall \rf{3.41} to infer that
$$
\rho_{\bf a} (B,{\fk M}) = ({\fk S}_{\bf a}-{\fk S}_{\bf a}(Q))
(J_{\bf a}(B) + E_{\bf a}) + O(Q^5 B^{s-1-k}A^{-1}),
$$
and hence that
$$
\rho_{\bf a} (B,{\fk M}) - {\fk S}_{\bf a} J_{\bf a}(B) \ll
|{\fk S}_{\bf a}(Q)| 
|J_{\bf a}(B) + E_{\bf a}| + {\fk S}_{\bf a} |E_{\bf a}| + Q^5
B^{s-1-k} A^{-1}.
$$
On the left hand side, we may invoke \rf{3.7}. On the right hand side,
we observe that by \rf{4.7} and \rf{3.5}, one has $J_a(B) + E_{\bf a}
\ll (a_1 \ldots a_s)^{-1/s} B^{s-k}$. Hence, we may  sum over ${\bf
  a}$ and apply Lemma 3.8, provided only that $s\ge k+2$, as we now
assume. Then
\be{4.8}
\sum_{|{\bf a}|\le A} |\rho_{\bf a} (B,{\fk M}) - {\fk S}_{\bf a}
J_{\bf a} B^{s-k} | \ll
A^{s-1} B^{s-k} Q^{-1/(2k)} + \Sigma + Q^5 B^{s-k-1} A^{s-1}.
\ee
where
\be{4.9}
\Sigma= \sum_{|{\bf a}| \le A} {\fk S}_{\bf a} |E_{\bf a}|.
\ee
Now, by \rf{4.7} followed by an application of H\"older's inequality,
$$
|E_{\bf a}| \le \int_{|\beta|\ge Q/(AB^k)}
|v_{\bf a} (\beta )| d\beta \le
\prod_{j=1}^s \Big(\int_{|\beta|\ge Q/(AB^k)} |v(a_j\beta)|^s
d\beta\Big)^{1/s}.
$$
However, whenever $0<|a|\le A$, then by \rf{3.2},
\begin{eqnarray*}
\int_{Q/(AB^k)}^\inf |v(a\beta)|^s d\beta &=&
\f{1}{|a|} \int_{Q|a|/(AB^k)}^\inf |v(\beta)|^s d\beta \\
&\ll&
|a|^{-1} B^{s-k} (1+Q|a| A^{-1})^{-s/k},
\end{eqnarray*}
which produces the estimate
$$
|E_{\bf a}|\ll |a_1 a_2 \ldots a_s|^{-1/s} B^{s-k}
\prod_{j=1}^s (1+Q |a_j| A^{-1})^{-1/k}.
$$
Hence, provided only that $A\ge Q$, Lemma 3.8 combined with 
a dyadic dissection argument for $|{\bf a}|$,
shows that
$$
\Sigma \ll Q^{-1/k} A^{s-1} B^{s-k}.
$$
We finally choose $\eta = \f{1}{6}$, and then by \rf{4.8}, conclude as
follows.
\\[2ex]
{\sc Lemma 4.2.} {\it Let $A\ge 1, B\ge 1$ and $Q=B^{1/6}$. Then,
  whenever $s\ge k+2$, one has}
$$
\sum_{|{\bf a}|\le A} |\rho_{\bf a}(B,{\fk M}) -{\fk S}_{\bf a} J_{\bf
  a} B^{s-k}| \ll A^{s-1} B^{s-k-1/(12k)}.
$$
\\[2ex]
{\bf 4.3. The minor arcs.} We begin the endgame with a variant of
Weyl's inequality.
\\[2ex]
{\sc Lemma 4.3.} {\it Let $A\ge 1, B\ge 1$, and suppose 
that $r\in\ZZ$ and $q\in \NN$ are
  coprime with $| \al - (r/q)|\le q^{-2}$. Then}
$$
\sum_{0< |a|\le A} |f(a\al)|^{2^{k-1}} \ll
AB^{2^{k-1}} \Big(\f{1}{q} + \f{1}{B} + \f{q}{AB^k}\Big) (AB q)^\eps.
$$
This is well known, but we give a brief sketch for completeness. Write
$K=2^{k-1}$. Then, as an intermediate step towards the ordinary form
of Weyl's inequality, one has
$$
|f(\beta)|^K \ll B^{K-1} + B^{K-k+\eps} \sum_{1\le h\le 2^k k!
  B^{k-1}} \min (B, \| h\beta \|^{-1})
$$
where $\|\beta \|$ denotes the distance of $\beta$ to the nearest
integer; compare the arguments underpinning Lemma 2.4 of Vaughan
\cite{hlm}. Now choose $\beta = a\al$ and sum over $a$. A divisor
function argument then yields
$$
\sum_{0<|a|\le A} |f(a\al)|^K \ll AB^{K-1} + B^{K-k} (AB)^\eps 
\sum_{h\ll AB^{k-1}} \min (B, \| h \beta \|^{-1}),
$$
and Lemma 4.3 follows from Lemma 2.2 of Vaughan \cite{hlm}.
\medskip

 Now let $\al \in{\fk m}$. By Dirichlet's theorem on diophantine
 approximations, there are $r\in\ZZ$, $q\in\NN$ with $q\le Q^{-1}
 AB^k$ and
$$
|q\al -r|\le Q (AB^k)^{-1}.
$$
But $\al \notin {\fk M}$, whence $q>Q$. Lemma 4.3 in conjunction with
H\"older's inequality now yields
\be{4.10}
\sup_{\al\in{\fk m}} \sum_{0<|a|\le A} |f(a\al)| \ll (AB)^{1+\eps}
Q^{-2^{1-k}}.
\ee
We now apply this estimate to establish the following.
\\[2ex]
{\sc Lemma 4.4.} {\it Let $s\in\NN$, $s= 2t + u$ with $t\in\NN$, $u=1$
  or $2$. Then there is a number $\del >0$ such that whenever $1\le
  B^k \le A\le B^{t-k}$ holds, then}
$$
\sum_{|{\bf a}|\le A} |\rho_{\bf a} (B,{\fk m}) |\ll A^{s-1}
B^{s-k-\del}.
$$
{\it Proof.} By \rf{4.5}, one has
$$
\sum_{|{\bf a}|\le A} |\rho_{\bf a}(B,{\fk m})| \le \int_{\fk m}
\Big( \sum_{0<|a|\le A} |f(a\al)|\Big)^s d\al.
$$
Moreover, by Cauchy's inequality and orthogonality,
\begin{eqnarray*}
\int_0^1 \Big( \sum_{0<|a|\le A} |f(a\al)|\Big)^{2t} d\al 
&\le &
(2A +1)^t \int_0^1 \Big(\sum_{0<|a|\le A} |f(a\al)|^2\Big)^t d\al \\
&\le &
(2A +1)^t U_t (A,B).
\end{eqnarray*}
On combining the last two inequalities with \rf{4.10} and Theorem 2.5,
we deduce that
\be{4.11}
\sum_{|{\bf a}| \le A} |\rho_{\bf a}(B,{\fk m})| \ll A^s B^{t+u-\del}
+ A^{s-1} B^{s-k-\del}
\ee
where any $0< \del < \f{1}{6} 2^{1-k}$ is admissible. Note that the
condition that $B^k \le A$ is required in Theorem 2.5, whereas the
inequality $A\le B^{t-k}$ makes the second term on the right of
\rf{4.11} the dominating one. This establishes the lemma.
\medskip

Theorem 1.1 is also available: one has ${\hat s}=2t$, and the theorem follows   on combining \rf{4.6} with Lemma 4.2
and Lemma 4.4.

\noindent
J\"org Br\"udern,\\ Mathematisches Institut, \\ Bunsenstrasse 3--5,\\ D-37073
G\"ottingen, Germany\\
bruedern@uni-math.gwdg.de
\\[1ex]
Rainer Dietmann,\\ Royal Holloway, University of London,\\
 Egham, Surrey TW20 0EX, 
United Kingdom\\
rainer.dietmann@rhul.ac.uk
\end{document}